\documentclass{article}
\usepackage{natbib}
\bibpunct{(}{)}{;}{a}{}{;} 
\title{Hybrid Pluralism}
\author{Andrew Aberdein\thanks{
Humanities and Communication,
Florida Institute of Technology,
150 West University Blvd, 
Melbourne, Florida 32901-6975, U.S.A.
{\sf aberdein@fit.edu}. }}
\renewcommand{\theequation}{$P_{\arabic{equation}}$}

\begin{document}

\maketitle

\begin{abstract}
Intuitively, a pluralist solution is one in which a single question receives multiple answers.  Such pluralist solutions have been proposed in many widely disparate contexts.  This paper restates the concept of pluralism with greater precision; distinguishes it from, and establishes its independence of, some other notions with which it is frequently confused; and briefly lays out some of the benefits that this more nuanced approach to pluralism may yield for the debates in which it may be invoked. 
\end{abstract}

Let $\mathcal{Q}$ be a set of questions: \[\mathcal{Q} = \{Q_1, Q_2, \dots, Q_n\}\] Let $\mathcal{A}$ be a set of answers: \[\mathcal{A} = \{A_1, A_2, \dots, A_m\}\] Let a \textit{position}, $P$, be a function with domain $\mathcal{Q}$ and codomain $\mathcal{A}$.  If $P$ is defined over $\mathcal{Q}$ then every question has a unique answer.  Where $P(Q_{i}) = A_{j}$ we will write $Q_{i}A_{j}$, read as ``$A_{j}$ answers $Q_{i}$''.  Thus we may write \[P = \{Q_{1}P(Q_{1}), Q_{2}P(Q_{2}), \dots, Q_{n}P(Q_{n})\}\]  Let $\mathcal{P}$ be a set of admissible positions, that is definitions of $P$.  Hence
\begin{eqnarray*}
\mathcal{P} = \{\{Q_{1}P_{1}(Q_{1}), Q_{2}P_{1}(Q_{2}), \dots, Q_{n}P_{1}(Q_{n})\}, \\
\{Q_{1}P_{2}(Q_{1}), Q_{2}P_{2}(Q_{2}), \dots, Q_{n}P_{2}(Q_{n})\},\\
\vdots \hspace{6.5em}\\
\{Q_{1}P_{k}(Q_{1}), Q_{2}P_{k}(Q_{2}), \dots, Q_{n}P_{k}(Q_{n})\}\}
\end{eqnarray*}  
where $1 \leq k \leq n^{m}$.  The elements of $\mathcal{P}$ may be simply ordered in an intuitive fashion.\footnote{Each element may be given the unique index
$1 + \mathop{\sum}_{i = 1}^{n}m^{n-i}(j_{i}-1)$, where $P(Q_{i}) = A_{j_{i}}$.}
If a pair of questions, $Q_{i}$, $Q_{j}$, receives the same answer in every element of $\mathcal{P}$ we may replace them in $\mathcal{Q}$ with their conjunction, $Q_{i} \wedge Q_{j}$, without loss of generality.\footnote{Alternatively, we may borrow from erotetic logic the convention that a question is defined by the answers it may receive, making $Q_{i}$, $Q_{j}$ and $Q_{i} \wedge Q_{j}$ equivalent if $Q_{i}$ and $Q_{j}$ receive the same answer in every position.}
We shall see that different constraints on the membership of $\mathcal{P}$ yield different conceptions of pluralism.

An asymmetry in this model is that multiple questions can receive the same answer, but multiple answers cannot hold for the same question (in the same position).  This is a consequence of the requirement that $P$ be a function.  This restriction may be removed by making $P$ a set-valued function, or multifunction, ranging over the power set of $\mathcal{A}$.  We might then distinguish positions in which multiple answers are assigned to single questions as those satisfying the condition
\begin{eqnarray*}
({\rm S})&(\exists x)(\exists y)(|P_{x}(Q_{y})| > 1)
\end{eqnarray*}
However, this is not pluralism: in pluralist situations a question may receive \textit{any} of a range of answers \textit{individually}; here a question would receive \textit{all} of a range of answers \textit{jointly}.  We shall rather call situations meeting (S) \textit{syncretist}.  In what follows, this possibility will be ignored, and $P$ will be assumed to be a function.  However, this constraint is not essential to any of the forthcoming definitions, and it will always remain possible to substitute the multifunction definition of $P$, permitting the presentation of syncretist variants of all the positions set out below.

The simplest cases arise when either $\mathcal{Q}$ or $\mathcal{A}$ is a singleton.  If $\mathcal{A}$ is the singleton set $\{A\}$ then 
\[\mathcal{P} = \{\{Q_{1}A, Q_{2}A, \dots, Q_{n}A\}\}\]
Since $\mathcal{P}$ has only one element, it is trivial that each $Q_{i}$ always has the same answer and thus that 
\[\mathcal{P} \cong \{\{(\mathop{\wedge}_{i = 1}^{n}Q_{i})A\}\}\]
Hence the case with singleton $\mathcal{A}$ is a special case of that with both $\mathcal{A}$ and $\mathcal{Q}$ singletons, $\mathcal{P} = \{\{QA\}\}$.  We shall say that this situation is \textit{global}, since in every position every question receives the same answer, and that it is \textit{monist}, since every question receives the same answer in every position.  More formally, we may say that $\mathcal{P}$ satisfies the conditions
\begin{eqnarray*}
({\rm G})&(\forall x)(\forall y)(\forall z)(P_{x}(Q_{y}) = P_{x}(Q_{z}))\\
({\rm M})&(\forall y)(\forall w)(\forall x)(P_{w}(Q_{y}) = P_{x}(Q_{y}))
\end{eqnarray*}
We shall call any situation in which (G) and (M) hold \textit{global monism}.   If $\mathcal{Q}$ is a singleton, $\{Q\}$, but $\mathcal{A}$ has more than one member then 
\[\mathcal{P} = \{\{QA_{1}\}, \{QA_{2}\}, \dots, \{QA_{m}\}\}\]
This situation is still global, but it is also \textit{pluralist}, since the same question receives different answers in different positions.  That is $\mathcal{P}$ satisfies the conditions
\begin{eqnarray*}
({\rm G})& (\forall x)(\forall y)(\forall z)(P_{x}(Q_{y}) = P_{x}(Q_{z}))\\
({\sim}{\rm M})&(\exists y)(\exists w)(\exists x)(P_{w}(Q_{y}) \neq P_{x}(Q_{y}))
\end{eqnarray*}
We shall call any situation meeting these conditions \textit{global pluralism}.

Now consider situations in which neither $\mathcal{A}$ nor $\mathcal{Q}$ are singletons.  If $\mathcal{P}$ is nevertheless constrained to a single element in which at least some of the questions receive different answers, say  \[P = \{Q_{1}A_{j_{1}}, Q_{2}A_{j_{2}}, \dots, Q_{n}A_{j_{n}}\}\] where $1 \leq j_{i} \leq m$, a new sort of situation will arise.  This situation, \textit{local monism}, satisfies the conditions
\begin{eqnarray*}
({\sim}{\rm G})&(\exists x)(\exists y)(\exists z)(P_{x}(Q_{y}) \neq  P_{x}(Q_{z}))\\
({\rm M})&(\forall y)(\forall w)(\forall x)(P_{w}(Q_{y}) = P_{x}(Q_{y}))
\end{eqnarray*}
(Where $m \geq n$ every question may receive a different answer; otherwise, two or more questions may be conjoined as described above.)    However, in many situations, including the maximal case where $\mathcal{P}$ includes every possible position $P$, there will be at least one  question which does not receive the same answer in every position, and at least one position will assign different answers to different questions.  That is to say, the conditions 
\begin{eqnarray*}
({\sim}{\rm G})&(\exists x)(\exists y)(\exists z)(P_{x}(Q_{y}) \neq  P_{x}(Q_{z}))\\
({\sim}{\rm M})&(\exists y)(\exists w)(\exists x)(P_{w}(Q_{y}) \neq P_{x}(Q_{y}))
\end{eqnarray*}
will be satisfied.  Call any situation meeting these conditions \textit{local pluralism}.  
It is now a simple exercise to demonstrate that every possible value of $\mathcal{P}$ corresponds to one of these four situations.

However, two additional constraints may be arrived at by weakening the conditions (G) and (M):
\begin{eqnarray*}
({\rm G^{\prime}})&(\exists x)(\forall y)(\forall z)(P_{x}(Q_{y}) = P_{x}(Q_{z}))\\
({\rm M^{\prime}})&(\exists y)(\forall w)(\forall x)(P_{w}(Q_{y}) = P_{x}(Q_{y}))
\end{eqnarray*}
\textit{Weak globalism} (G$^{\prime}$) states that there are \textit{some} positions in which every question receives the same answer and \textit{weak monism} (M$^{\prime}$) states that there are \textit{some} questions which receive the same answer in every position.  Conversely, the negations of (G$^{\prime}$) and (M$^{\prime}$) represent strengthenings of ($\sim$G) and ($\sim$M):
\begin{eqnarray*}
({\sim}{\rm G}^{\prime})&(\forall x)(\exists y)(\exists z)(P_{x}(Q_{y}) \neq  P_{x}(Q_{z}))\\
({\sim}{\rm M}^{\prime})&(\forall y)(\exists w)(\exists x)(P_{w}(Q_{y}) \neq P_{x}(Q_{y}))
\end{eqnarray*}
\textit{Strong localism} $(\sim$G$^{\prime})$ states that every position contains questions which receive different answers; \textit{strong pluralism} $(\sim$M$^{\prime})$ states that no question receives the same answer in every position.
Since the sets $\mathcal{A}$ and $\mathcal{Q}$ are non-empty, it can be easily verified that (G) $\Rightarrow $ (G$^{\prime}$), (M) $\Rightarrow $ (M$^{\prime}$), ($\sim$G) $\wedge$ (M) $\Rightarrow $ ($\sim$G$^{\prime}$) and (G) $\wedge$ ($\sim$M) $\Rightarrow $ ($\sim$M$^{\prime}$).  Hence, these additional constraints only have an effect when both ($\sim$G) and ($\sim$M) apply, that is in cases of local pluralism.  Constraints on $\mathcal{P}$ such that at least one question, but not every question, receives the same answer in every position meet the conditions ($\sim$G), ($\sim$M) and (M$^{\prime}$).  We will call this situation \textit{hybrid pluralism}.  Conversely, the conditions ($\sim$G), ($\sim$M) and (G$^{\prime}$) are satisfied when at least one, but not every, position assigns the same answer to every question.  We will call this situation \textit{hybrid localism}.  (Note that some situations will exhibit both hybrid forms.)  Finally, we may refer to situations in which  ($\sim$G), ($\sim$M), ($\sim$G$^{\prime}$) and ($\sim$M$^{\prime}$) all apply as \textit{strict local pluralism}.

The foregoing may be made clearer by a simple example.  Let $\mathcal{Q} = \{Q_1, Q_2, Q_3\}$ and $\mathcal{A} = \{A_1, A_2\}$.  Then, maximally,
\renewcommand{\theequation}{$P_{\arabic{equation}}$}
\begin{eqnarray}
\mathcal{P} = \{\{Q_{1}A_{1}, Q_{2}A_{1}, Q_{3}A_{1}\}, \\
\{Q_{1}A_{1}, Q_{2}A_{1}, Q_{3}A_{2}\}, \\
\{Q_{1}A_{1}, Q_{2}A_{2}, Q_{3}A_{1}\},\\
\{Q_{1}A_{1}, Q_{2}A_{2}, Q_{3}A_{2}\},\\
\{Q_{1}A_{2}, Q_{2}A_{1}, Q_{3}A_{1}\},\\
\{Q_{1}A_{2}, Q_{2}A_{1}, Q_{3}A_{2}\},\\
\{Q_{1}A_{2}, Q_{2}A_{2}, Q_{3}A_{1}\},\\
\{Q_{1}A_{2}, Q_{2}A_{2}, Q_{3}A_{2}\}\}
\end{eqnarray}    
Restricting $\mathcal{P}$ to either $\{P_{1}\}$ or $\{P_{8}\}$ produces global monism.  Including both elements, that is $\mathcal{P} = \{P_{1}, P_{8}\}$, produces global pluralism.  Each of $\{P_{2}\}$ through $\{P_{7}\}$ produces local monism.  Any combination which includes at least one of both $P_{1}$ and $P_{8}$, and $P_{2}$ through $P_{7}$ will be hybrid localist.  Combinations of elements, such as $\{P_{1}, P_{2}\}$ or $\{P_{5}, P_{6}, P_{7}, P_{8}\}$, in which at least one question always receives the same answer, represent hybrid pluralism.  All other sets of elements, including the unrestricted $\mathcal{P}$, represent strict local pluralism.

It should now be obvious that the global/local distinction is not identical to that between monism and pluralism.  This point has been recognized in previous taxonomies of pluralism, such as that of Susan Haack \citeyearpar{Haack78} and Michael Resnik  \citeyearpar{Resnik96}.  However, for Haack, localism and globalism are subdivisions of pluralism, whereas we have seen that the two distinctions are independent.  This idiosyncrasy of Haack's classification leads to her misleadingly describing local monism as local \textit{pluralism} \citep[p.~223]{Haack78}. Neither she nor \citet[p.~499]{Resnik96}, who adopts her definition, considers the position which I call local pluralism.

Conventional presentations of pluralism offer a stark choice between (implicitly global) monism and  (implicitly global) pluralism.  Even the more sophisticated Haack/Resnick taxonomy does little to alleviate this dichotomy, since its other distinction is subordinated to that between monism and pluralism.
Hence many debates in which a pluralist approach is mooted may be characterized as a slide from global monism to global pluralism.
Global monism represents an ideal scenario, in which a consensus has been achieved on a common solution to an all-encompassing problem.  Unfortunately, this ideal is remote.  This failure of global monism then leads inevitably to a confrontation between two or more broadly defined and frequently irreconcilable positions, that is global pluralism.

Local and hybrid pluralism represent attractive intermediate stages, at which the slide towards global pluralism may be arrested.  Local pluralism captures the intuitive thought that many large debates subsume smaller debates.  Hybrid pluralism adds that consensus may be found in some of the smaller debates, even amongst proponents of sharply contrasted  positions in the larger debate.  In the remainder of the paper I will discuss examples of this more nuanced approach.

Examples from recent ethical thought include Ronald Dworkin's \citeyearpar{Dworkin} analysis of the abortion debate and Martha Nussbaum's \citeyearpar{NUS00} discussion of the interaction of feminism and multiculturalism.  Dworkin argues that if the abortion debate is characterized crudely, for example as the question of whether the foetus is a person, then the two sides will be irreconcilable.  But, if the question is broken down into smaller sub-questions, then an unexpected degree of consensus may be found: most members of both lobbies agree that the lives of both the mother and the foetus have value, and in some circumstances they agree on whether abortion should be permissible. Nussbaum responds in a similar way to the apparent conflict between feminism and the multiculturalist admonition to respect non-Western cultures---despite their anti-feminist practices.  She breaks down the global pluralism of crude multiculturalism by challenging the monolithic conception of `culture'.  This permits a local pluralism, which Nussbaum pushes towards hybrid pluralism by arguing for consensus over some key intercultural questions.  Whether or not one accepts Dworkin or Nussbaum's analyses, they both exhibit the utility of hybrid pluralism.

Another area in which pluralism has been especially conspicuous in recent years is the philosophical discussion of the proliferation of non-classical logics.  Most recent accounts of logical pluralism, such as \citep{Restall00}, are implicitly globalist: the pluralist offers multiple answers to the single question ``What is the true logic?''.
However, most important non-classical systems have a substantial common subsystem: classical logic itself.  Intuitionists accept the classical logicality of decidable propositions; paraconsistentists, consistent propositions; quantum logicians, compatible propositions; and so forth.   It has even been suggested that this `recapture' of classical logic is a necessary criterion of logicality.\footnote{`Perhaps \dots any genuine `logical system' should contain classical logic as a special case' \citep[p.~135]{vBE94}.}  In the light of these concessions, it is more plausible to characterize logical pluralism as an instance of hybrid pluralism, with the proponents of the different systems agreeing over the formalization of a central common discourse, and diverging in other areas.

Even classical logic and its conservative extensions, that is the system(s) of logic characteristically endorsed by logical monists, may be seen to exhibit a degree of hybrid pluralism. 
Many classicists say that classical logic is the `one true logic'.  The natural understanding of this remark is global monism, with the unique answer understood as (first-order) classical logic. However, only a minority of classicists would defend a restriction of logicality to first-order classical logic \citep{QUI53, HAZ99}. Most recognize a variety of quantified or modal conservative extensions as equally logical. Taking this intuition seriously, while retaining global monism, would require the single formal system to somehow combine all 
the extensions of classical logic for which there might ever be a need. Yet despite some na\"{\i}vely misplaced optimism, the construction of such a compound system is a task of formidable technical difficulty if more than a small range of familiar extensions are to be used  (see, for example, \citealt{GAB98}). Furthermore, most conceivable applications would employ extensions containing only some of the extra material rather than the unwieldy compound system containing it all. So local monism would seem a closer approximation to the actual commitments of the classical programme (\textit{cf.\ }\citealt[p.~44]{HAA74}). 
However, it is widely acknowledged that many discourses addressed by conservative extensions of classical logic, including alethic and other modalities, lack an unambiguous choice of formalization. This suggests monism about the discourse of first-order classical logic, and some of its extensions, and pluralism about that of some other extensions: that is to say, hybrid pluralism.

\end{document}